\documentclass[12pt,a4paper,reqno]{amsart}
\usepackage{amssymb}
\usepackage{multirow}
\usepackage{graphicx}
\usepackage{amsmath,amsbsy,amsfonts,amsthm,latexsym,amsopn,amstext,amsxtra,euscript,amscd,mathrsfs}

\theoremstyle{plain}

\theoremstyle{definition}

\theoremstyle{remark}

\begin{document}

\title[Generalized Cosecant Numbers and the Hurwitz Zeta Function]{Generalized Cosecant Numbers and the Hurwitz Zeta Function}


\author{Victor Kowalenko}
\address{Department of Mathematics and Statistics,
The University of Melbourne, Victoria 3010, Australia}
\email{vkowa@unimelb.edu.au}

%    General info
\subjclass[2010]{11M06,11M35,11P83,11Y55,11Z05}

\date{\today}

\keywords{coefficient, generalized cosecant number, multiplicity, partition, Pochhammmer notation, symmetric 
polynomial, zeta function}

\date{\today}

\begin{abstract}
This paper presents recent developments concerning the generalized cosecant numbers $c_{\rho,k}$, which emerge as the 
coefficients of the power series expansion for the important fundamental function $z^{\rho}/\sin^{\rho} z$. These 
coefficients can be computed for all, including complex, values of $\rho$ by using the relatively novel graphical 
method known as the partition method for a power series expansion or by using intrinsic routines in Mathematica. In
fact, they represent polynomials in $\rho$ of degree $k$, where $k$ is the power of $z$. In addition, though related 
to the Bernoulli numbers, they possess more properties and do not diverge like the former. The partition method for
a power series expansion has the advantage that it yields the $k$-behaviour of the highest order coefficients. Thus, general 
formulas for such coefficients are derived by considering the properties of the highest part partitions. It is then shown 
how the generalized cosecant numbers are related to the specific symmetric polynomials that arise from summing over 
quadratic powers of integers. Consequently, integral values of the Hurwitz zeta function for even powers are expressed for 
the first time ever in terms of ratios of the generalized cosecant numbers.    
\end{abstract}

\maketitle

\section{Introduction}
The cosecant numbers, $c_k$, are defined in Ref.\ \cite{kow11} as the rational coefficients of the power series 
expansion for cosecant, in particular by
\begin{eqnarray} 
\csc z= \frac{1}{\sin z} \equiv \sum_{k=0}^{\infty} c_k \, z^{2k-1} \;.
\label{one}\end{eqnarray}
Via the discrete graphical method known as partition method for a power series expansion a general formula for them is derived in terms 
of all the integer partitions summing to $k$, which is   
\begin{equation}
c_{k}= (-1)^k \sum_{\substack{\lambda_1,\lambda_2,\lambda_3,\dots,\lambda_k
\\\sum_{i=1}^k i\, \lambda_i=k}}^{k,\lfloor k/2 \rfloor, \lfloor k/3
\rfloor,\dots,1} (-1)^{N_k} N_k!  \prod_{i=1}^k \Bigl(
\frac{1}{(2i+1)!} \Bigr)^{\lambda_i}  \frac{1}{\lambda_i!}\, .
\label{two}\end{equation}
In this result $\lambda_i$ represents the multiplicity or the number of occurrences of each part $i$ in the partitions, while the sum 
of the multiplicities or length of the partition is represented by $N_k$, i.e.\ $N_k = \sum_{i=1}^k \lambda_i$. For the partitions 
summing to $k$, the multiplicity of a part $i$ ranges from zero to $\lfloor k/i \rfloor$, where $\lfloor x \rfloor$ denotes
the floor function or the greatest integer less than or equal to $x$. Generally, when $k$ is large, most of the multiplicities for 
a partition vanish as we shall see shortly. Using \eqref{two} one finds that $c_0=1$, $c_1=1/6$, $c_2=7/360$, $c_4=31/3\cdot 7!$ 
and so on. In addition, the reader should observe that an equivalence symbol has been introduced into \eqref{one} because
the power series expansion on the rhs is divergent for $|z| \geq \pi$, while the lhs is always defined. For these values the 
rhs must be regularized in the manner described in Refs.\ \cite{kow001}-\cite{kow11b}. Nevertheless, since the power series is convergent 
for the other values of $z$, it is permissible to replace the equivalence symbol by an equals sign. Then one obtains the standard
power series for $x \csc x$, which is given as No.\ 1.411(11) in Ref.\ \cite{gra94}. Consequently, the cosecant numbers can be expressed
in terms of the Bernoulli numbers as 
\begin{eqnarray}
c_k=\frac{(-1)^{k+1}}{(2k)!}\,(2^{2k}-2)\,B_{2k} \;\;.
\label{twoa}\end{eqnarray}
As discussed in Ref.\ \cite{kow11}, the cosecant numbers possess far more properties than their more famous counterparts. They have also
the major advantage that they do not diverge like them. In fact, in most situations such as \ref{tw0a} and the Euler Maclaurin summation 
formula, the divergence of the Bernoulli numbers is tamed by the factor of $(2k)!$ in the denominator. Thus, one is often 
required to divide two very large numbers by each other, which is avoided when using the cosecant numbers. Hence the above
equation can be regarded as awkward or even clumsy. Nevertheless, for the interested reader it should be mentioned that Ref.\ \cite{kow09a}  
shows how the Bernoulli numbers and polynomials can be evaluated by using the partition method for a power series expansion.  

A better method of expressing the cosecant numbers is in terms of the Riemann zeta function. By using No.\ 
9.616 in Ref.\ \cite{gra94}, we find that
\begin{eqnarray}
c_k=2 \bigl( 1-2^{1-2k} \bigr) \, \frac{\zeta(2k)}{\pi^{2k}} \;,
\label{three}\end{eqnarray}
where $\zeta(2k)$ represents the Riemann zeta function. Thus, we observe that the $c_k \approx 2 \pi^{-2k}$ or $\pi^{2k} 
\approx 2/c_k$ for $k \gg 1$, although this approximation gives the misleading impression that the $c_k$ are irrational. 
Consequently, \eqref{two} becomes another method of determining even integer values of the Riemann zeta function. 

Ref.\ \cite{kow11} not only presents numerous applications of cosecant numbers, but also demonstrates how they are related
to other numbers such as the secant numbers and, more importantly, how the sets of the resulting numbers can be generalized or 
extended by introducing an arbitrary power $\rho$ to their generating function. Specifically, the generalized cosecant numbers are 
given by
\begin{eqnarray} 
\csc^{\rho} z = \frac{1}{\sin^{\rho} z} \equiv \sum_{k=0}^{\infty} c_{\rho,k} \, z^{2k-\rho} \;,
\label{four}\end{eqnarray}
where
\begin{equation}
c_{\rho,k}= (-1)^k \sum_{\substack{\lambda_1,\lambda_2,\lambda_3,\dots,\lambda_k
\\\sum_{i=1}^k i\, \lambda_i=k}}^{k,\lfloor k/2 \rfloor, \lfloor k/3
\rfloor,\dots,1} (-1)^{N_k} (\rho)_{N_k} \prod_{i=1}^k \Bigl(
\frac{1}{(2i+1)!} \Bigr)^{\lambda_i}  \frac{1}{\lambda_i!}\, .
\label{five}\end{equation}
In \eqref{five} $(\rho)_{N_k}$ denotes the Pochhammer notation for $\Gamma(\rho+N_k)/\Gamma(\rho)$, where $\Gamma(x)$ represents
the gamma function.

\section{Computation}
To calculate the generalized cosecant numbers via \eqref{five}, we need to determine the specific contribution made by 
each integer partition that sums to $k$. For example, if we wish to evaluate $c_{\rho,6}$, then we require all the contributions 
made from the partitions summing to 6, which appear in the first column of Table\ \ref{table1}. Each part in a partition 
is assigned a specific value, which depends on the function being studied. In the case of $x^{\rho}/\sin^{\rho}x$, the part
$i$ is assigned a value of $(-1)^{i+1}/(2i+1)!$. In addition, since each part occurs $\lambda_i$ times in a partition, we need 
to multiply $\lambda_i$ values or calculate $(-1)^{(i+1)\lambda_i}/((2i+1)!)^{\lambda_i}$. The second column in Table\ \ref{table1} 
displays the multiplicities of the parts in all the partitions summing to 6, while the third column presents the length $N_k$ 
for each partition. Thus, we see that most of the multiplicities vanish as stated earlier.

\begin{table}
%\begin{center}
\small
\begin{tabular}{|c|c c c c c c|c|} \hline
Partition & $\lambda_1$ & $\lambda_2$ & $\lambda_3$ & $\lambda_4$ & $\lambda_5$ & $\lambda_6$ & $N_k$     \\ \hline
$\{6\}$ & & & & & &  1 & 1\\
$\{5,1\}$ &1 & & & & 1 & & 2    \\
$\{4,2\}$ &  &1 & & 1 & & & 2 \\
$\{4,1,1\}$ &2 &  & & 1 & & & 3 \\
$\{3,3\}$ & &  & 2 & &  & & 2 \\
$\{3,2,1\}$ & 1 & 1 & 1 & & & & 3 \\
$\{3,1,1,1\}$ & 3 & & 1 & & & & 4 \\
$\{2,2,2\}$ &  & 3 & & &  & & 3 \\
$\{2,2,1,1\}$ & 2 & 2& & & & & 4 \\
$\{2,1,1,1,1\}$ & 4 & 1 & & & & & 5\\
$\{1,1,1,1,1,1\}$ & 6 & & & & & & 6\\
\hline
\end{tabular}
%\end{center}
%\normalsize
\vspace{0.3cm}
\caption{Multiplicities of the partitions summing to 6}
\label{table1}\end{table}

Associated with each partition is a multinomial factor that is determined by taking the factorial of $N_k$ and dividing by the 
factorials of all the multiplicities. E.g., for the partition $\{2,1,1,1,1\}$ in Table\
\ref{table1}, we have $\lambda_1=4$ and $\lambda_2=1$, while the other multiplicities vanish. Hence, the multinomial factor  
becomes $5!/(4!\, 1!)=5$. When the function is accompanied by an arbitrary power, say $\rho$, a further modification must 
be made. Each partition is then multiplied by the Pochhammer factor of $\Gamma(N_k+\rho)/\Gamma(\rho)$ divided by $N_k!$. That is,
for each partition we must include the extra factor of $(\rho)_{N_k}/N_k!$. For $\rho=1$, this simply yields unity and thus, the 
multinomial factor remains unaffected. Consequently, \eqref{five} reduces to \eqref{two} for $\rho=1$. For $\rho=2$ we obtain the 
cosecant-squared numbers as given in Theorem\ 5 of Ref.\ \cite{kow11}, which are given by
\begin{align}
c_{2,k} & = (-1)^k \sum_{\substack{\lambda_1,\lambda_2,\lambda_3,\dots,\lambda_k
\\\sum_{i=1}^k i\, \lambda_i=k}}^{k,\lfloor k/2 \rfloor, \lfloor k/3
\rfloor,\dots,1} (-1)^{N_k} (2)_{N_k} \prod_{i=1}^k \Bigl(
\frac{1}{(2i+1)!} \Bigr)^{\lambda_i}  \frac{1}{\lambda_i!}
\nonumber\\
& = \;\; \frac{(2k-1)}{\bigl( 1- 2^{1-2k} \bigr)}\; c_k \;\;.
\label{fivea}\end{align}
In the above result one can replace $(2)_{N_k}$ by $(N_k+1) (1)_{N_k}$. Hence we obtain \eqref{two} 
again except $N_k !$ is now multiplied by $N_k+1$. On the other hand, if $\rho=-1$, then the coefficients 
are simply equal to the power series expansion for $\sin z$ divided by $z$. Therefore, we arrive at
\begin{equation}
\sum_{\substack{\lambda_1,\lambda_2,\lambda_3,\dots,\lambda_k
\\\sum_{i=1}^k i\, \lambda_i=k}}^{k,\lfloor k/2 \rfloor, \lfloor k/3
\rfloor,\dots,1} (-1)^{N_k} (\rho)_{N_k} \prod_{i=1}^k \Bigl(
\frac{1}{(2i+1)!} \Bigr)^{\lambda_i}  \frac{1}{\lambda_i!} = \frac{1}{(2k+1)!} \, .
\label{fiveb}\end{equation}
Hence, we have an expression for the reciprocal of $(2k+1)!$ in terms of a sum over partitions summing to $k$.

If we examine \eqref{five} more closely, then we see that the product deals with calculating the contribution made by 
each partition based on the values of the multiplicities, while the sum refers to all partitions summing to $k$. Hence the sum 
covers the range of values for each multiplicity. For example, $\lambda_1$ attains a maximum value of $k$, which corresponds to 
the partition with $k$ ones, while $\lambda_2$ attains a maximum value of $[k/2]$, which corresponds to the partition with 
$[k/2]$ twos in it. For odd values of $k$, the partition with $[k/2]$ twos also possesses a one, i.e. $\lambda_1=1$. Thus, 
it can be seen that the maximum value of $\lambda_i$ is always $[k/i]$, which becomes the upper limit for each multiplicity 
in both \eqref{two} and \eqref{five}. Moreover, each partition in the sums must satisfy the constraint, $\sum_{i=1}^k i\lambda_i=k$.

As an example, let us calculate $c_{\rho,6}$. According to Table\ \ref{table1} there are eleven partitions summing to 6.
Therefore, we need to determine eleven contributions in the sum over the partitions. By applying the steps mentioned above 
to \eqref{five}, we find that the contributions from the partitions in the same order as the table are
\begin{align}
 c_{\rho,6}& = -(\rho)_1 \;\frac{1}{13!} + \frac{(\rho)_2}{2!}\;
\frac{2!}{1! \cdot 1!}\; \frac{1}{3! \cdot 11!} +\frac{(\rho)_2}{2!} \;
\frac{2!}{1! \cdot 1!}\; \frac{1}{5! \cdot 9!} - \frac{(\rho)_3}{3!}\;
\frac{3!}{1! \cdot 2!} \; \frac{1}{3!^2 \cdot 9!}
\nonumber\\
&  +\;\; \frac{(\rho)_2}{2!}\; \frac{2!}{2!}\; \frac{1}{(7!)^2}-
\frac{(\rho)_3}{3!}\; \frac{3!}{1! \cdot 1! \cdot 1!}\; \frac{1}{3! \cdot 5! \cdot 7!} +
\frac{(\rho)_4}{ 4!}\, \frac{4!}{1! \cdot 3!} \; \frac{1}{(3!)^3 \cdot 7!} -
\frac{(\rho)_3}{3!}\; \frac{3!}{3!} \;\frac{1}{(5!)^3} 
\nonumber\\
& + \;\; \frac{(\rho)_4}{ 4!}\, \frac{4!}{2! \cdot 2!} \; \frac{1}{(3!)^2 \cdot (5!)^2}
- \frac{(\rho)_5}{ 5!}\, \frac{5!}{4! \cdot 1!} \;\frac{1}{(3!)^4 \cdot 5!}
+ \frac{(\rho)_6}{ 6!}\, \frac{6!}{6!} \;\frac{1}{(3!)^6} \;\;.
\label{six}\end{align}
The interesting property of the above result is that when the length $N_k$ appearing in the Pochhammer terms is even, the contribution 
from the partition is positive while if it is odd, then the contribution is negative. This behaviour applies to all even
values of $k$. On the other hand, if $k$ is odd, then the contributions with an odd number of parts are positive, while those
from an even number of parts are negative. Furthermore, by introducing \eqref{six} into Mathematica \cite{wol92} and wrapping it entirely 
around the combination of the Simplify and Expand routines, one obtains
\begin{align}
c_{\rho,6} & = \frac{1}{5884534656000}\, \Bigl(1061376 \rho + 3327584 \rho^2 + 4252248 \rho^3 + 2862860 \rho^4 
\nonumber\\
& + \;\; 1051050 \rho^5 + 175175 \rho^6 \Bigr) \;.
\label{seven}\end{align}
Since the denominator equals $2/(9 \cdot 15!)$, we arrive at the $k=6$ result in Table\ \ref{table2}.

\begin{table}
\begin{center}
\small
\begin{tabular}{|c|l|} \hline
$k$ & $c_{\rho,k}$    \\ \hline
$0$ &  $ 1 $  \\ [0.1 cm] 
$1$   &  $\frac{1}{3!} \; \rho $ \\ [0.1 cm]
$2$  &  $ \frac{2}{6!} \; \bigl(2\rho+5\rho^2 \bigr)$  \\ [0.1 cm]
$3$ &  $ \frac{8} {9!}\;\bigl( 16 \rho+42 \rho^2+35 \rho^3 \bigr) $ \\ [0.1 cm]
$4$ &  $ \frac{2}{3 \cdot 10!} \; \bigl( 144 \rho+404\rho^2+420 \rho^3+175 \rho^4 \bigr)$ \\ [0.1 cm]
$5$ & $ \frac{4}{3 \cdot 12!} \; \bigl( 768 \rho+2288 \rho^2+2684 \rho^3+1540 \rho^4+ 385 \rho^5 \bigr)$ \\ [0.1 cm]
$6$ & $ \frac{2}{9 \cdot 15!} \; \bigl( 1061376 \rho+3327594 \rho^2+4252248 \rho^3+2862860 \rho^4+1051050 
\rho^5  +175175 \rho^6 \bigr)$ \\ [0.1 cm]
$7$  & $ \frac{1} {27 \cdot 15!} \; \bigl( 552960 \rho+ 1810176 \rho^2+ 2471456 \rho^3 + 1849848 \rho^4+ 820820 \rho^5 + 210210 \rho^6$ \\
& $ +25025 \rho^7 \bigr)$ \\ [0.1 cm]
$8$ & $ \frac{2}{45 \cdot 18!} \bigl( 200005632 \rho + 679395072 \rho^2+ 978649472 \rho^3 + 792548432 \rho^4 + 397517120 \rho^5 $ \\ 
& $+ 125925800 \rho^6 + 23823800 \rho^7+ 2127125 \rho^8 \bigr)$ \\ [0.1 cm]
$9$ & $ \frac{4}{81 \cdot 21!} \bigl(129369047040 \rho+ 453757851648 \rho^2 + 683526873856 \rho^3 +  589153364352 \rho^4 $ \\
& $ + 323159810064 \rho^5 + 117327450240 \rho^6 + 27973905960 \rho^7 + 4073869800 \rho^8 +  282907625 \rho^9 \bigr)$ \\[0.1 cm]
$10$ & $\frac{2}{6075 \cdot 22!} \bigl( 38930128699392 \rho + 140441050828800 \rho^2 + 219792161825280 \rho^3 + 199416835425280 \rho^4 $ \\
& $ + 117302530691808 \rho^5 + 47005085727600 \rho^6 + 12995644662000 \rho^7 + 2422012593000 \rho^8$ \\ 
& $ + 280078548750 \rho^9 + 15559919375 \rho^{10} \bigr)$  \\[0.1 cm]
$11$ & $ \frac{8}{243 \cdot 25!} \bigl(494848416153600 \rho + 1830317979303936 \rho^2 + 2961137042841600 \rho^3 $ \\
& $ + 2805729689044480 \rho^4 + 1747214980192000 \rho^5 + 755817391389984 \rho^6 + 232489541684400 \rho^7$ \\
& $ + 50749166067600 \rho^8 + 7607466867000 \rho^9 + 715756291250 \rho^{10} + 32534376875 \rho^{11} \bigr)$ \\[0.1 cm]
$12$ & $ \frac{2}{2835 \cdot 27!} \bigl( 1505662706987827200 \rho + 5695207005856038912 \rho^2 + 9487372599204065280 \rho^3$ \\
& $ + 9332354263294766080 \rho^4 + 6096633539052376320 \rho^5 + 2806128331871953088 \rho^6 $ \\
& $ + 937291839756592320 \rho^7 + 229239926321406000 \rho^8 + 40598842049766000 \rho^9 $ \\
& $ + 5005999501002500 \rho^{10} + 390802935022500 \rho^{11} + 14803141478125 \rho^{12} \bigr) $ \\[0.1 cm]
$13$ & $ \frac{232}{81 \cdot 30!} \bigl( 844922884529848320 \rho + 3261358271400247296 \rho^2 +5576528334428209152 \rho^3$ \\
& $+5668465199488266240 \rho^4 + 3858582205451484160 \rho^5 + 1870620248833400064 \rho^6 $ \\
& $+ 667822651436228288 \rho^7 + 178292330746770240 \rho^8 + 35600276746834800 \rho^9 $ \\
& $+ 5225593531158000 \rho^{10} + 539680243602500 \rho^{11} + 35527539547500 \rho^{12} + 1138703190625 \rho^{13} \bigr) $ \\[0.1 cm]
$14$ & $ \frac{2}{1215 \cdot 30!}\; \bigl(138319015041155727360 \rho + 543855095595477762048 \rho^2 + 952027796641042464768 \rho^3 $ \\
& $+ 996352286992030556160 \rho^4 + 703040965960031795200 \rho^5 + 356312537387839432192 \rho^6 $ \\
& $+ 134466795172062184832 \rho^7 + 38526945410311117760 \rho^8 + 8436987713444690400 \rho^9 $\\
& $+ 1404048942958662000 \rho^{10} + 173777038440005000 \rho^{11} + 15258232341852500 \rho^{12} $\\
& $+ 858582205731250 \rho^{13} +23587423234375 \rho^{14} \bigr)$ \\ [0.1 cm]
$15$ & $\frac{1088}{729 \cdot 35!} \; \bigl( 562009739464769840087040 \rho + 2247511941596311764074496 \rho^2 $ \\
& $+ 4019108379306905439830016 \rho^3 + 4317745925208072594259968 \rho^4 $ \\
& $+ 3145163776677939429416960 \rho^5 +1656917203539032341530624 \rho^6 $ \\
& $+ 655643919364420586023424 \rho^7 + 199227919419039256217472 \rho^8 $ \\
& $+ 46995751664475880185920 \rho^9 + 8614026107092938211680 \rho^{10} +1214778349162323946000 \rho^{11} $ \\
& $+ 128587452922193265000 \rho^{12} + 9720180867524627500 \rho^{13} + 472946705787806250 \rho^{14}$ \\
& $+ 11260635852090625 \rho^{15} \bigr) $ \\[0.1 cm]
\hline
\end{tabular}
\end{center}
%\end{center}
\normalsize
\vspace{0.3cm}
\caption{Generalized cosecant numbers $c_{\rho,k}$ up to $k=15$}
\label{table2}\end{table}

Table\ \ref{table2} displays the generalized cosecant numbers up to $k=15$, which have been obtained by introducing the multiplicities 
of all partitions summing to $k$ into the sum in \eqref{two}. For $k\!>\!10$, the partition method for a power series expansion becomes 
laborious due to the exponential increase in the number of partitions. To circumvent this problem, a general computing methodology has been
developed in Refs.\ \cite{kow17} and \cite{kow12}, which is based on representing all the partitions summing to a specific order $k$ by a 
partition tree and invoking the bivariate recursive central partition (BRCP) algorithm. From this computing methodology general 
expressions for the coefficients of any power series expansion obtained. For example, the symbolic form in the case of the partitions
summing to 6 is given by \vspace{0.25 cm} \newline 
DS[6]:= p[6] q[1] a + p[1] p[5] q[2] a$^{\wedge}$(2) 2! + p[1]$^{\wedge}$(2) p[4] q[3] a$^{\wedge}$(3) 3!/2! + \newline
p[1]$^{\wedge}$(3) p[3] q[4] a$^{\wedge}$(4) 4!/3! + p[1]$^{\wedge}$(4) p[2] q[5] a$^{\wedge}$(5) 5!/4! + p[1]$^{\wedge}$(6) q[6] 
a$^{\wedge}$(6) + \newline
p[1]$^{\wedge}$(2) p[2]$^{\wedge}$(2) q[4] a$^{\wedge}$(4) 4!/(2! 2!) + p[1] p[2] p[3] q[3] a$^{\wedge}$(3) 3! + \newline + 
p[2] p[4] q[2] a$^{\wedge}$(2) 2! + p[2]$^{\wedge}$(3) q[3] a$^{\wedge}$(3) + p[3]$^{\wedge}$(2) q[2] a$^{\wedge}$(2) $\quad .$
\vspace{0.25 cm}\newline
Such an expression can easily be imported into Mathematica \cite{wol92}. Then the coefficients of the inner series p[k] are set equal to
the assigned values of the parts. For the generalized cosecant numbers, this means we set \vspace{0.25 cm} \newline
p[k$_{-}$]:=(-1)$^{\wedge}$(k+1)/(2k+1)! $\;\;,$ \vspace{0.25 cm}\newline
while the q[k], which are referred to as the coefficients of the outer series, are set equal to the coefficients of the binomial series, 
viz. \vspace{0.25 cm} \newline
q[k$_{-}$]:=Pochhammer[$\rho$,k]/k! $\;\;.$ \vspace{0.25 cm}\newline
In addition, we set the parameter a equal to unity, thereby obtaining $c_{\rho,6}$. On the other hand, changing
p[k] to \vspace{0.25 cm} \newline 
p[k$_{-}$] := (-1)$^{\wedge}$(k + 1)/(2 k)! $\;\;,$ \vspace{0.25 cm} \newline
will yield an entirely different set numbers known as the generalized secant numbers $d_{\rho, k}$ \cite{kow11,kow17}. These numbers represent
the coefficients in the power series expansion for $\sec^{\rho} z$. In terms of the partition method for a power series expansion they
are given by
\begin{eqnarray}
d_{\rho,k}= (-1)^k \sum_{\substack{ n_1, n_2,n_3, \dots,n_k=0 \\ \sum_{i=1}^k i n_i=k}}^{k,[k/2],[k/3],\dots,1}   
(-1)^N \, (\rho)_N \prod_{i=1}^k \left(\frac{1}{(2i)!} \right)^{\!n_i} \frac{1}{n_i !} \;\;.
\label{sevene}\end{eqnarray}
Thus, the symbolic representation for DS[6] is not only general, but also powerful.

It should also be mentioned that Mathematica \cite{wol92} via its SeriesCoefficient routine is able to determine the generalized cosecant
numbers very quickly due to the fact that it has been optimized. For example, to obtain the first fifteen generalized cosecant numbers, 
one need only type:\vspace{0.2cm}\newline
\vspace{0.2cm}
In[1]:= Table[SeriesCoefficient[(z/Sin[z])$^{\wedge}\rho$, \{z, 0, k\}], \{k, 0, 30, 2\}].\newline
The result for $k=7$ generated by this command is:\vspace{0.2 cm} \newline
Out[2]= (8191 $\rho$)/37362124800 + (15019 (-1 + $\rho$) $\rho$)/12454041600 + (517457 (-2 + $\rho$) (-1 + $\rho$) $\rho$)/301771008000 + (
 169 (-3 + $\rho$) (-2 + $\rho$) (-1 + $\rho$) $\rho$)/183708000 + (83 (-4 + $\rho$) (-3 + $\rho$) (-2 +$\rho$) (-1 + $\rho$) $\rho$)/391910400 
+ ( 7 (-5 + $\rho$) (-4 + $\rho$) (-3 + $\rho$) (-2 + $\rho$) (-1 + $\rho$) $\rho$)/335923200 + ((-6 + $\rho$) (-5 + $\rho$) (-4 + $\rho$) 
(-3 + $\rho$) (-2 + $\rho$) (-1 + $\rho$) $\rho$)/1410877440 .
\vspace{0.2 cm}\newline
Although there are far less contributions than by the partition method for a power series expansion, this form for 
the generalized cosecant numbers is still cumbersome because one must, once again, employ the Expand and Simplify routines in 
Mathematica to obtain the forms displayed in Table\ \ref{table2}. Moreover, by using this approach one cannot determine the 
$k$-dependence of the coefficients of the polynomials, whereas the partition method for a power series expansion is able to reveal 
this behaviour for the highest order coefficients as explained in the next section.

\section{Coefficients of the Generalized Cosecant Numbers}
Despite the fact that the contributions in \eqref{five} alternate in sign according to whether the length $N_k$ for each partition is even 
or odd, the final forms for the generalized cosecant numbers only possess positive coefficients. Moreover, the highest order terms in the 
$c_{\rho,k}$ is $O(\rho^k)$. For example, from \eqref{seven} we see that $c_{\rho,6}$ is a sixth order polynomial in $\rho$. 
The term that is responsible for the $\rho^6$ term in \eqref{five} emanates from the partition with six ones in it because it yields $(\rho)_6$. 
Consequently, we see that the generalized cosecant numbers are polynomials of degree $k$ with fixed coefficients. We can express them 
as $c_{\rho,k}= \sum_{i=1}^k C_{k,i} \rho^i$. Our aim in this section is to determine the highest order coefficients as a function of $k$.

Since we know that the highest order term in the generalized cosecant numbers is determined by the partition with the most number of ones, 
i.e. $k$ ones, we evaluate the contribution from this partition in \eqref{five}, which is $(\rho)_k/(3!)^k k!$. Hence the highest order 
term in $\rho$ is the coefficient of $\rho^k$ in the Pochhammer factor, which is
\begin{eqnarray}
C_{k,k}= \frac{1}{(3!)^k k!} \;\;.
\label{seven1}\end{eqnarray}

The second highest order term, namely $C_{k,k-1}$, is the sum of the contributions from two partitions. First, there is
the $\rho^{k-1}$ term from the partition with $k$ ones and second, there is the highest order term from the partition with
$k-2$ ones and one two. The first term represents the coefficient of $\rho^{k-1}$ in the Pochhammer factor in the previous
calculation, while the second term represents the $\rho^{k-1}$ power in $-(\rho)_{k-1}/(5! \cdot (3!)^{k-1} (k-2)!)$.
Combining the contributions yields
\begin{eqnarray}
C_{k,k-1}= \frac{1}{(3!)^k k!} \sum_{i=1}^{k-1}i -\frac{1}{ 5\cdot (3!)^{k-2} (k-2)!}= \frac{1}{5 \cdot
(3!)^k  (k-2)!} \;\;.
\label{seven2}\end{eqnarray}
The next highest order term or $C_{k,k-2}$ is the sum of the contributions from four partitions, viz.\ the $\rho^{k-2}$ power
from the partition with $k$ ones, the second leading order term from the partition with $k-2$ ones and one two and the
leading order terms in the partitions with $k-3$ ones and one three and $k-4$ ones and two twos. 

To evaluate each of these contributions in general form, we require the formula that gives the coefficient of $\rho$ to an 
arbitrary power in each Pochhammer factor. According to Chapters\ 24 and 18 of Refs.\ \cite{abr65} and \cite{spa87} 
re\-spec\-tive\-ly, the Poch\-hammer polynomials can be expressed as
\begin{eqnarray}
(y)_k= \frac{\Gamma(y+k)}{\Gamma(y)}=(-1)^k \sum_{j=0}^{k} (-1)^j s^{(j)}_k \,y^j \;\;,
\label{seven3}\end{eqnarray}
where the integers $s^{(j)}_k$ are known as the (signed) Stirling numbers of the first kind with
$s^{(0)}_{k} \!=\! s^{(k)}_0 \!=\!0$ for $k \!\geq \!1$ and  $s_0^0 \!=\!1$. It should be mentioned
that the Stirling numbers of the first kind are often represented as $s(k,j)$. They also satisfy the
following recurrence relation:
\begin{eqnarray}
s^{(j)}_{k+1}=s^{(j-1)}_k-k \, s^{(j)}_{k} \;\;.
\label{seven4}\end{eqnarray}
In the appendix of Ref.\ \cite{kow09a}, a general formula for these numbers is derived, which is given as
\begin{eqnarray}
s_k^{(k-j)}= (-1)^j \sum_{i_j=j}^{k-1} i_j \sum_{i_{j-1}=j-1}^{i_j-1} \sum_{i_{j-2}=j-2}^{i_{j-1}-1} i_{j-2}
\cdots \sum_{i_1=1}^{i_2-1}i_1 \;\;. 
\label{seven5}\end{eqnarray}
By using this result we can calculate general formulas for specific values of $j$ ranging from 0 to 4. These are
\begin{eqnarray}
s^{(k)}_k=1 \quad, \quad s^{(k-1)}_k= -\binom{k}{2}\quad, \quad s^{(k-2)}_k = \frac{(3k-1)}{4}\, \binom{k}{3} \;\;,
\nonumber\\
s^{(k-3)}_{k} = - \binom{k}{2} \binom{k}{4} \quad , \quad s^{(k-4)}_k = \frac{1}{48} \bigl(15k^3-30 k^2 +5k +2 \bigr)
\binom{k}{5} \;\;.
\label{seven6}\end{eqnarray}
Note that the denominator of the last result has been corrected here compared with Ref.\ \cite{kow09a}, where 336 appears
instead of 48. These results have been obtained by introducing the nested sum formula given by (\ref{seven5})
as a module in Mathematica \cite{wol92}. Thus, in general one obtains a polynomial in $k$ of degree $2j$ with one of the terms
in the polynomials being equal to $(-1)^j \binom{k}{j+1}$. In Chapter\ 6 of Ref.\ \cite{kow17} the partition method
for a power series expansion is used to derive an alternative formula for the Stirling numbers of the first kind, which will 
enables general formulas for $s^{(k-\ell)}_k$ to be determined more easily than using (\ref{seven5}). Then it is found that 
the Stirling numbers of the first kind can be represented as $s^{(k-\ell)}_k = (-1)^{\ell} \binom{k}{\ell+1} r_{\ell}(k)$, where 
the $r_{\ell}(k)$ are polynomials of degree $\ell-1$ and are displayed in Table\ \ref{table4}. For odd values of $\ell$ the 
polynomials possess a common external factor of $k(k-1)$. As an aside, it should be mentioned that the 
corresponding results for the Stirling numbers of the second kind are also derived in Ref.\ \cite{kow17}. 

\begin{table}
%\begin{center} 
{\renewcommand{\arraystretch}{1.2}
\centering{ \fontsize{12 pt}{14 pt}\selectfont
\begin{tabular}{|c|l|} \hline
$\ell$ &  $ r_{\ell}(k)$  \\ \hline
$1$ &  $1$  \\
$2$ & $ \frac{1}{4} (3k-1) $ \\[0.1 cm]
$3$ & $ \frac{1}{2}\, k (k-1) $  \\ [0.1 cm]
$4$ & $\frac{1}{48} \bigl( 15 k^3-30k^2+5 k+2 \bigr) $ \\ [0.1 cm]
$5$ & $ \frac{1}{16} \, k(k-1) \bigl( 3k^2-7k -2 \bigr) $ \\ [0.1 cm]
$6$ & $\frac{1}{576} \bigl(63 k^5 -315 k^4 +315 k^3 +91 k^2 -42 k -16 \bigr) $ \\ [0.1 cm]
$7$ & $\frac{1}{144} \, k (k-1) \bigl( 9k^4- 54 k^3+ 51 k^2 + 58 k + 16  \bigr) $ \\ [0.1 cm]
$8$ & $\frac{1}{3840} \bigl( 135 k^7 -1260 k^6 +3150 k^5 -840 k^4 -2345 k^3 -540 k^2  -66262636 k$ \\
& $+ 596367504  \bigr) $ \\ [0.1 cm]
$9$ & $ \frac{1}{768} \, k (k-1) \bigl( 15 k^7 -180 k^6 +630 k^5 -448 k^4 -665 k^3 +100 k^2 +404 k -144 \bigr) $ \\ [0.1cm]
$10$ & $\frac{1}{9216} \bigl( 99 k^9 -1485 k^8 +6930 k^7 -8778 k^6 -8085 k^5 + 8195 k^4 +11792 k^3 $ \\
& $ +2068 k^2 -2288 k -768 \bigr)$
\\ [0.1 cm] \hline
\end{tabular}

}
}
%\end{center}
\caption{The polynomials $r_{\ell}(k)$ in the Stirling numbers of the first kind}
\label{table4}
\end{table}

With the above results $C_{k,k-2}$ reduces to
\begin{align}
C_{k,k-2} & = \frac{s^{(k-2)}_k}{(3!)^k k!} + \frac{s^{(k-2)}_{k-1}}{5! \cdot (3!)^{k-2} (k-2)!}+ 
\frac{s^{(k-2)}_{k-2}}{7! \cdot (3!)^{k-3} (k-3)!} 
\nonumber\\
& + \;\; \frac{s^{(k-2)}_{k-2}}{ 2!\cdot (5!)^2  (3!)^{k-4} (k-4)!} \;\;. 
\label{seven6a}\end{align}
Introducing the results in \eqref{seven6}, we find that the coefficients reduce to
\begin{eqnarray}
C_{k,k-2}= \frac{21 k+17}{175 \, (3!)^{k+1} \,(k-3)!} \;\;.
\label{seven7}\end{eqnarray}

The expressions for $C_{k,k-\ell}$ become more difficult to evaluate as $\ell$ increases. However, from the above results,
we see that a pattern is developing. First, the power of $3!$ appears to be increasing as the power of $\rho$ decreases. 
Next the factorial in the denominator decrements by unity. That is, $C_{k,k-1}$ goes as $1/ (3!)^k (k-2)!$, while $C_{k,k-2}$ 
goes as $1/(3!)^{k+1} (k-3)!$. Furthermore, the numerator for $C_{k,k-1}$ is constant, whereas it is linear in $k$ for
$C_{k,k-2}$. So we conjecture that the next coefficient is given by
\begin{eqnarray}
C_{k,k-3}= \frac{ak^2 +bk+c}{(3!)^{k+2} (k-4)!} \;\;.
\label{seven8}\end{eqnarray}
The reason for choosing $1/(k-\ell-1)!$ in the denominator for $C_{k,k-\ell}$ is that the first value of these coefficients
only begins when $k= \ell+1$. Thus, we can put $k=4$ in the above result and equate it to $C_{4,1}$ in Table\ \ref{table2},
which equals 144/5443200. Similarly, we put $k=5$ and $k=6$ in (\ref{seven8}) and equate the values respectively to $C_{5,2}$
and $C_{6,3}$ in the same table. Then we arrive at the following set of equations:
\begin{gather}
16 a+ 4b+c =\frac{216}{175}\;\;,\\
25a+5b+c=\frac{312}{175}\;\;,\\
36a+6b+c=\frac{2124}{875} \;\;.
\label{seven9}\end{gather}
The solution to the above set of equations is $a=6/125$, $b=102/875$ and $c=0$, which means in turn that $C_{k,k-3}$ is 
given by
\begin{eqnarray}
C_{k,k-3}= \frac{k^2 +17k/7}{125 (3!)^{k+1} (k-4)!} \;\;.
\label{seven10}\end{eqnarray}
Putting $k=8$ in this formula yields 73/26453952000, which agrees with the coefficient of 397517120 in Table\
\ref{table2} for $k=8$ when it is multiplied by the external factor of $2/45 \cdot 18!$. In using this approach it did not 
matter whether our conjecture had the correct power of $3!$ in the denominator initially provided the dependence upon $k$ 
was correct. However, it is crucial that the correct factorial appears in the denominator. 

In the case of $C_{k,k-4}$ it would be conjectured by
\begin{eqnarray}
C_{k,k-4}= \frac{ak^3 +bk^2 +ck +d}{(3!)^{k+4} (k-5)!} \;\;.
\label{seven11}\end{eqnarray}
Then a set of four linear equations is required with $k$ ranging from 5 to 8. Solving the equations using
the LinearSolve routine in Mathematica \cite{wol92} yields
\begin{eqnarray}
C_{k,k-4}= \frac{3 \bigl(3k^3 +102 k^2/7+ 289 k/49 -11170/539\bigr) }{625\, (3!)^{k+3} \, (k-5)!} \;\;.
\label{seven12}\end{eqnarray}
Finally, putting $k=9$ into the above result yields a value of 229051/733303549440000, which agrees with the coefficient of
$\rho^5$ for $k=9$ in Table\ \ref{table2}.

Determining general formulas for the coefficients of the lowest order terms in the generalized cosecant numbers is a much
more formidable problem. If we consider the preceding methods, then to obtain the equivalent of (\ref{seven2}) for
the lowest order coefficient, viz. $C_{k,1}$, we need to evaluate the contributions from all partitions whereas previously we
only needed a fixed number, e.g., four for determining $C_{k,k-2}$. This is clearly not possible since the number of partitions is 
not fixed, but increases exponentially. In addition, the conjectural approach breaks down when $k$ appears in the second subscript 
of the coefficients. For example, to determine $C_{k,2}$ via this approach, $\ell$ would now be equal to $k-2$. 

For $|\rho| \gg k$ we may use the first four coefficients derived above as a means of approximating the generalized cosecant numbers. 
That is, the generalized cosecant numbers can be approximated by 
\begin{align}
c_{\rho, k} \overset{|\rho| \gg k}{\approx} \frac{\rho^k}{(3!)^k k!} + \frac{\rho^{k-1}}{5 \cdot (3!)^k (k-2)!} + 
\frac{(21 k+17) \rho^{k-2}}{175 \cdot (3!)^{k+1} (k-3)!} + \frac{(k^2+17k/7) \rho^{k-3}}{125 \cdot (3!)^{k+1}(k-4)!}\;.
\label{sevenf}\end{align} 
To gain an appreciation of this approximation, let us denote the ratio of \eqref{sevenf} to the corresponding values of $c_{\rho,k}$ in Table\
\ref{table2} by $\beta(\rho,k)$. Table\ \ref{table3} presents values of $\beta(\rho,k)$ for integer values of $\rho$ ranging from 10 to 1000. 
They have been given to six decimal places with no rounding-off. From the table we see that when $\rho$ is close to $k$ or smaller, which occurs
towards the right hand top corner, $\beta(\rho,k)$ is not close to unity, but for all other values, it is. Therefore, provided $\rho$ is 
significantly greater than $k$, \eqref{sevenf} represents an accurate approximation for the generalized cosecant numbers.   

\begin{table}
\begin{center}
\small
\begin{tabular}{|c|l|l|l|l|l|} \hline
$\rho$ & $k=6$ & $k=8$ & $k=10$ & $k=12$ & $k=15$    \\ \hline
$10$ &  $ 0.998905$ & $0.985655$ & $0.929830$ & $0.801477$ & $0.501086$  \\ [0.1 cm] 
$15$   &  $ 0.999741$ & $0.996144$ & $0.977941$ & $0.925497$ & $0.751944$ \\ [0.1 cm]
$20$  &  $0.999910 $ & $0.998571$ & $0.991119$ & $0.966957$ & $ 0.870459$   \\ [0.1 cm]
$30$ &  $0.999981 $ & $0.999669$ & $0.997755$ & $ 0.990752$ & $0.956671$  \\ [0.1 cm]
$50$ & $0.999997 $ & $0.999951$ & $0.999644$ & $0.998405$ & $ 0.991292$ \\ [0.1 cm]
$100$ & $ 0.999999$ & $ 0.999997$ & $0.999974$ & $0.999676$ & $ 0.992370$ \\ [0.1 cm] 
$1000$ & $0.999999 $ & $ 0.999999$ & $ 0.999999$ & $0.999999$ & $0.999999$ \\ [0.1 cm]
\hline
\end{tabular}
\end{center}
%\end{center}
\normalsize
\vspace{0.3cm}
\caption{The ratio $\beta(\rho,k)$ for various values of $\rho$ and $k$}
\label{table3}\end{table}

There is one interesting feature in the table that needs to be mentioned. If one examines the ratios when (1) $k=10$ and $\rho=20$
and (2) $k=15$ and $\rho=30$, then it is readily observed that the ratio is more accurate for the first case than in the second case despite
the fact they are both good approximations. This means that as $k$ increases, $|\rho/k|$ must also increase in order to obtain the same
value of $\beta(\rho,k)$ for the lower values of $k$. For example, the value of $\beta(20,10)$ is about 0.991, which in the case of $k=15$
only is only reached when $\rho$ is about 50. Hence in the $k=10$ case it only takes twice the value of $k$ to achieve the same value
of $\beta(\rho,k)$ as the $k=15$ case, which requires at least three times the value of $k$. This has ramifications when the 
relationship between the generalized cosecant numbers and the Hurwitz zeta function is discussed next. 

\section{Connection to Hurwitz zeta function}
It is found in Refs. \cite{fon17} and \cite{fon16}  that the generalized cosecant numbers can also be expressed as
\begin{equation}
c_{2v,i} = 2^{2i}\, \frac{ \Gamma(2v-2i)}{\Gamma(2v)} \, s(v,i)  \;,\quad i<v\;,
\label{nine}\end{equation}
where $s(v,n)$ represents the $n$th elementary symmetric polynomial obtained by summing over quadratic powers or squared integers, viz.\ 
$1^2,2^2,\ldots,(v-1)^2$. That is,
\begin{equation}
s(v,n)=\sum_{1 \leq i_1<i_2<\cdots<i_n<v-1} x_{i_1} x_{i_2} \cdots
x_{i_n}\;\;,
\label{ten}\end{equation}
where $x_{i_1} <x_{i_2}< \cdots <x_{i_n}$ and each $x_{i_j}$ is equal to at least one value in the set 
$\left\{1,2^2,3^2,\dots,(v-1)^2 \right\}$. For the three lowest values of $n$ the symmetric polynomials are
$$
s(v,0)=1 \;,  \quad s(v,1)=(v-1)v(2v-1)/6 \; , $$
and
\begin{equation} 
s(v,2)= \frac{(5 v+1)}{4 \cdot 6!} \; (2v-4)_5 \; ,
\label{eleven}\end{equation}
while for the four highest values of $n$ they are given by
$$s(v,v-1) = (v-1)!^2 \;, \quad
s(v,v-2)=  (v-1)!^2\,\bigl( \zeta(2) -\zeta(2,v) \bigr) \;,$$
$$ s(v,v-3) = \frac{(v-1)!^2}{2} \Bigl( (\zeta(2) -\zeta(2,v))^2+
\zeta(4,v)-\zeta(4) \Bigr) \;, $$
and
\begin{align}
s(v,v-4) & = \frac{(v-1)!^2}{6} \Bigl( (\zeta(2) -\zeta(2,v))^3-3
(\zeta(4)-\zeta(4,v))(\zeta(2) -\zeta(2,v)) 
\nonumber\\
& + \;\; 2 (\zeta(6) -\zeta(6,v))\Bigr)\; .
\label{twelve}\end{align}

The $n= v-\ell$ results have been obtained by beginning with the general form given by \eqref{ten}. For $s(v,v-4)$ this becomes
\begin{eqnarray}
s(v,v-4)=  \frac{1}{6} \prod_{i=1}^{v-1} i^2 \sum_{\substack{ j_1,j_2,j_3=1 \\ j_1\neq j_2 \neq j_3}}^{v-1} \frac{1}{j_1^2 j_2^2 j_3^2} \;.
\label{thirteen}\end{eqnarray} 
To evaluate this result, the constraint that none of the $j_i$ is equal to another must be removed. This is accomplished
by dropping it and eliminating all the possibilities where at least one of $j_i$ is equal to another. Consequently, we must 
consider subtracting all the cases where two of the indices are equal to another and finally when three indices are equal to one 
another. The product over $i$ yields $\Gamma(v)^2$. Thus, we arrive at
\begin{eqnarray}
s(v,v-4) = \frac{1}{2} \Gamma(v)^2\Bigl( \sum_{j_1,j_2,j_3=1}^{v-1} \frac{1}{j_1^2 j_2^2 j_3^2} -3 \sum_{j_1,j_2=1}^{v-1} 
\frac{1}{j_1^4 j_2^2}+2 \sum_{j_1=1}^{v-1} \frac{1}{j_1^6} \Bigr) \;\;,
\label{fourteen}\end{eqnarray}
which yields the result given by \eqref{twelve}. In a similar manner one finds that
\begin{align} 
s(v,v-5) & = \frac{1}{24} \; \Gamma(v)^2 \Bigl( \bigl( \zeta(2)-\zeta(2,v) \bigr)^4- 6 \bigl( \zeta(4) -\zeta(4,v) \bigr)
\bigl( \zeta(2) -\zeta(2,v) \bigr)^2
\nonumber\\
&  + \; 8 \bigl(\zeta(6) -\zeta(6,v) \bigr) \bigl( \zeta(2) -\zeta(2,v) \bigr) +3 \bigl( \zeta(4) -\zeta(4,v) \bigr)^2 \Bigr) 
\nonumber\\
& - \;6  \bigl(\zeta(8) -\zeta(8,v) \bigr) \Bigr)\;.
\label{fourteena}\end{align}
In general, $s(v,v-\ell)$ is given by
\begin{align}
s(v,v-\ell) =  \sum_{\substack{j_1, \dots,j_{\ell-1}=1 \\ j_1< j_2,\cdots<j_{\ell-1}}}^{v-1} \prod_{i=1}^{v-1} \prod_{k=1}^{\ell-1} \frac{i^2}{j_k^2}
= \frac{\Gamma(v)^2}{(\ell-1)!} \sum_{ \substack{j_1, \dots,j_{\ell-1}=1 \\ j_1 \neq j_2 \cdots j_{n-2} \neq j_{\ell-1}}}^{v-1}
 \prod_{k=1}^{\ell-1} \frac{1}{j_k^2}\;.
\label{fourteenb}\end{align}
In order to solve the sum on the rhs, one needs to remove the constraint that each of the $j_i$ cannot equal one another. This means that
we need to subtract all the possibilities when at least one of the indices is equal to one another from the sum where all the indices can equal 
each other other, viz.\ $\sum_{j_1,\dots,j_{\ell-1}=1}^{v-1} \prod_{k=1}^{\ell-1}1/j_k^2$. The number of sums that appear on the rhs becomes
$p(\ell-1)$, where $p(k)$ is the partition function or the number of partitions summing to $k$. E.g., in \eqref{fourteen} we see that there
are three sums for $s(v,v-4)$, since $p(3)$ is equal to three as a result of the partitions, \{1,1,1\}, \{2,1\} and \{3\}. In addition, the powers
of the $j_i$ will be twice the parts in the partitions. That is, the sum $\sum_{j_1,j_2=1}^{v-1} 1/(j_1^4 j_2^2)$ corresponds to the partition
\{2,1\}. Moreover, the sign outside each sum is $(-1)^{N_k +1}$, while the factors preceding each sum is the factorial of the length $N_k$ 
divided by the factorials of the multiplicities of each part and the parts taken to the power of their multiplicities. For example, in the
case of the partition \{2,1\}, the sum of the parts is three, while the multiplicities are both equal to 1. Hence the factor outside the sum 
becomes $3!/ (2\cdot 1! \cdot 1 \cdot 1!)=3$. See Ref.\ \cite{fon16} for more details.

If we introduce the results for $s(v,v-n)$ into \eqref{nine}, then for $n=1$ we obtain
\begin{equation}
\frac{\Gamma(v)}{\Gamma(v+1/2)} = \frac{2}{\sqrt{\pi}}\; c_{2v,v-1} \;\;,
\label{fifteen}\end{equation}
which is only valid for $v>1$. Alternatively, \eqref{fifteen} can be expressed as $B(v,1/2) = 2\, c_{2v, v-1}$, where $B(x,y)$
represents the beta function. Furthermore, according to No.\ 2.5.3.1 in Ref.\ \cite{pru86}, we have
\begin{eqnarray}
\int_0^{\pi/2} 
\left\{ \begin{matrix} \sin x \\ \cos x \end{matrix} \right\}^{2v-1}  
\; dx= \frac{\sqrt{\pi}}{2}\, 
\frac{\Gamma(v)}{\Gamma(v+1/2)} = c_{2v, v-1}\;\;,
\label{fifteena}\end{eqnarray} 
while from Appendix\ I.1.9 of the same reference, we find that
\begin{eqnarray}
\left\{ \begin{matrix} \sin x \\ \cos x \end{matrix} \right\}^{2v-1}  =\frac{1}{2^{2v-2}}
\sum_{k=0}^{v-1}(\mp 1)^{v-k-1} \binom{2v-1}{k} \left\{\begin{matrix} \sin (2v-2k-1)x \\ \cos (2v-2k-1)x 
\end{matrix} \right\} .
\label{fifteenb}\end{eqnarray}
A surprising property of the above analysis is that the $c_{2v,k}$ are coefficients in the power series expansion
of $(x/\sin x)^{2v}$, but in \eqref{fifteena} they are related to positive powers of sine and cosine taken to $2v-1$ when
$k= v-1$. Inserting the second result into the first one yields
\begin{eqnarray}
\frac{1}{2^{2v-2}}\sum_{k=0}^{v-1}\frac{(-1)^{v-k-1}}{(2v-2k-1)}\, \binom{2v-1}{k} =c_{2v,v-1}\; .
\label{fifteenc}\end{eqnarray}
This result can be checked by implementing it in Mathematica as follows:\vspace{0.2cm} \newline
C2vvminus1[v$_{-}$] := 2$^{\wedge}$(2 - 2 v) Sum[(-1)$^{\wedge}$(v - k - 1) Binomial[2 v - 1, k]/(2 v - 2 k - 1), \{k, 0, v - 1\}]
\vspace{0.2cm}\newline
Putting v equal to 5 in this line of code yields 128/315, while putting $\rho=10$ in the $k=4$ result of Table\ \ref{table2} gives the
same value.

It emerges that the terms at the upper limit of the sum in \eqref{fifteenc} contribute more to the value of $c_{2v,v-1}$ than those
at the lower limit. Therefore, we express \eqref{fifteenc} as
\begin{eqnarray}
c_{2v,v-1} =2^{2-2v}\, \binom{2v-1}{v}\sum_{k=0}^{v-1}\frac{(-1)^{k}}{(2k+1)}\, \prod_{j=1}^k \Bigl( \frac{1-j/v}{1+j/v} \Bigr)\; .
\label{fifteend}\end{eqnarray}
Now we expand the terms in the sum and product, thereby obtaining
\begin{align}
c_{2v,v-1}& =2^{2-2v}\, \binom{2v-1}{v} \Bigl[ \Bigl(1- \frac{2}{v}+ \frac{2}{v^2} +O \bigl(\frac{1}{v^3} \bigr) \Bigr)-
\frac{1}{3}  \Bigl(1- \frac{2}{v}+ \frac{2}{v^2} +O \bigl(\frac{1}{v^3} \bigr) \Bigr) \Bigl(1- \frac{4}{v} 
\nonumber\\
& + \;\;\frac{8}{v^2} +O \bigl(\frac{1}{v^3} \bigr) \Bigr)  + \frac{1}{5}  \Bigl(1- \frac{2}{v}+ \frac{2}{v^2} +O \bigl(\frac{1}{v^3} \bigr) 
\Bigr) \Bigl(1- \frac{4}{v}+ \frac{8}{v^2} +O \bigl(\frac{1}{v^3} \bigr) \Bigr)
\nonumber\\
& \times \;\;  \Bigl(1- \frac{6}{v}+ \frac{18}{v^2} +O \bigl(\frac{1}{v^3} \bigr) \Bigr) - \cdots \Bigr]\; .
\label{fifteene}\end{align}
In more compact notation the above result can be written as
\begin{eqnarray}
c_{2v,v-1}=2^{2-2v}\, \binom{2v-1}{v} \Bigl( \sum_{j=0}^{v-1} \frac{(-1)^j}{(2j+1)}- \frac{1}{v}\sum_{j=0}^{v-1} \frac{(-1)^j}{(2j+1)}\;
(j+1)(j+2) + O\bigl( \frac{1}{v^2}\bigr) \Bigr)\;. 
\label{fifteenf}\end{eqnarray}
The first sum is a known result given by No.\ 4.1.3.4 in Ref.\ \cite{pru86}, while the second sum requires decomposition. Then we arrive at
\begin{eqnarray}
\sum_{j=0}^{v-1} \frac{(-1^j}{2j+1}\; (j+1)(j+2)= \frac{1}{2} \sum_{j=0}^{v-1}(-1)^j j+ \frac{5}{4} \sum_{j=0}^{v-1}(-1)^j
+ \frac{3}{4} \sum_{j=0}^{v-1} \frac{(-1)^j}{2j+1} \;\;.
\label{fifteeng}\end{eqnarray}
The last sum in \eqref{fifteeng} is simply another occurrence of the first sum in \eqref{fifteenf}. Consequently, $c_{2v, v-1}$ becomes
\begin{align}
c_{2v,v-1}& =2^{2-2v}\, \binom{2v-1}{v} \Bigl( \frac{\pi}{4} + \frac{(-1)^{v-1}}{2} \, \beta(v+1/2)+ \frac{(-1)^{v-1}}{2v} \,\lfloor v/2 \rfloor
\nonumber\\
& -\;\;\frac{5}{8v}\bigl( 1- (-1)^v \bigr) + \frac{3}{4v} \frac{(-1)^{v-1}}{2} \, \beta(v+1/2)  + O\bigl( \frac{1}{v^2} \bigr)  \Bigr)\;, 
\label{fifteenh}\end{align}
where $\lfloor x \rfloor$ denotes the floor function or the greatest integer less than or equal to $x$, $\beta(v)= (\psi((v+1)/2)-\psi(v/2))/2$ 
and $\psi(x)$ is the digamma function. For large values of $v$, the leading term given by the $\pi/4$ term on the rhs is a good asymptotic 
approximation to $c_{2v,v-1}$.

For $n=2$, \eqref{nine} yields
\begin{equation}
\sum_{k=1}^{v-1} \frac{1}{k^2} = \frac{\pi^2}{6}- \zeta(2,v) = \frac{2}{3} \frac{c_{2v,v-2}}{c_{2v,v-1}} \;\;,
\label{sixteen}\end{equation}
where $\zeta(x,y)$ represents the Hurwitz zeta function and $v>2$. By adopting the same approach for $n=3$ to $n=6$, we 
arrive at
\begin{equation}
\sum_{k=1}^{v-1} \frac{1}{k^4} = \frac{\pi^4}{90}- \zeta(4,v) = \frac{4}{9} \Bigl(\frac{c_{2v,v-2}}{c_{2v,v-1}}\Bigr)^2 -
\frac{4}{15} \frac{c_{2v,v-3}}{c_{2v,v-1}}\;\;,
\label{seventeen}\end{equation}
\begin{equation}
\sum_{k=1}^{v-1} \frac{1}{k^6} = \frac{\pi^6}{945}- \zeta(6,v) = \frac{4}{105} \;\frac{c_{2v,v-4}}{c_{2v,v-1}} -
\frac{4}{15} \frac{c_{2v,v-3}}{c_{2v,v-1}} \; \frac{c_{2v,v-2}}{c_{2v,v-1}} +\frac{8}{27} \Bigl(\frac{c_{2v,v-2}}{c_{2v,v-1}}\Bigr)^3\;,
\label{eighteen}\end{equation}
\begin{align}
\sum_{k=1}^{v-1} \frac{1}{k^8} & = \frac{\pi^8}{9450}- \zeta(8,v) = \frac{8}{14175} \left( 350\, \Bigl(\frac{c_{2v,v-2}}{c_{2v,v-1}}\Bigr)^4 -
420 \,  \frac{c_{2v,v-3}}{c_{2v,v-1}} \; \Bigl(\frac{c_{2v,v-2}}{c_{2v,v-1}}\Bigr)^2 \right.
\nonumber\\
 &\left.  +\; 63 \Bigl(\frac{c_{2v,v-3}}{c_{2v,v-1}}\Bigr)^2 +\; 60 \, \frac{c_{2v,v-4}}{c_{2v,v-1}}  \, \frac{c_{2v,v-2}}{c_{2v,v-1}} 
-5  \, \frac{c_{2v,v-5}}{c_{2v,v-1}} \right)\;,
\label{nineteen}\end{align}
and
\begin{align}
\sum_{k=1}^{v-1} \frac{1}{k^{10}} & = \frac{\pi^{10}}{93555}- \zeta(10,v) = \frac{4}{93555} \left( 3080\, 
\Bigl(\frac{c_{2v,v-2}}{c_{2v,v-1}}\Bigr)^5 -
4620 \,  \frac{c_{2v,v-3}}{c_{2v,v-1}} \; \Bigl(\frac{c_{2v,v-2}}{c_{2v,v-1}} \Bigr)^3 \right.
\nonumber\\
 &  +\; 1386 \Bigl(\frac{c_{2v,v-3}}{c_{2v,v-1}}\Bigr)^2 \frac{c_{2v,v-2}}{c_{2v,v-1}} +\; 660 \, 
\frac{c_{2v,v-4}}{c_{2v,v-1}}  \, \Bigl(\frac{c_{2v,v-2}}{c_{2v,v-1}}\Bigr)^2 
-198  \, \frac{c_{2v,v-4}}{c_{2v,v-1}} \;\frac{c_{2v,v-3}}{c_{2v,v-1}} 
\nonumber\\
& \left. -\; 55 \, \frac{c_{2v,v-5}}{c_{2v,v-1}} \;\frac{c_{2v,v-2}}{c_{2v,v-1}} +3\, \frac{c_{2v,v-6}}{c_{2v,v-1}} \right) \;.
\label{twenty}\end{align}
where $v>2$ for the first result, $v>3$ in the second, etc. In principle, this process can be continued for higher powers of the sum
on the lhs by determining larger values of $\ell$ in the symmetric polynomials, $s(v,v-\ell)$. Consequently, we see that integer values 
of the Hurwitz zeta function for even powers can now be expressed in terms of ratios of the generalized cosecant numbers, which is 
indeed fascinating in view of the intractability of this famous function. Moreover, in the limit as $v\to \infty$, we obtain new results 
for the Riemann zeta function such as
\begin{equation}
\zeta(4) =\frac{\pi^4}{90}=\lim_{v \to \infty} \Bigl\{  \frac{4}{9} \Bigl(\frac{c_{2v,v-2}}{c_{2v,v-1}}\Bigr)^2 - \frac{4}{15} 
\frac{c_{2v,v-3}}{c_{2v,v-1}} \Bigr\} \;\;. 
\label{twentyone}\end{equation}
Unfortunately, we cannot introduce \eqref{sevenf} into the above result because $|\rho/k|$ is approximately equal to 2 when $k$ is equal 
to either $v-1$ and $v-2$, whereas we have observed in Table \ref{table3} that as $\rho$ or $2v$ increases, the ratio $2v/k$ needs
to increase dramatically in order to ensure that $\beta(\rho,k)$ remains close to unity. Otherwise, \eqref{sevenf} does not represent
a good approximation for the generalized cosecant numbers. Therefore, we require asymptotic forms as in \eqref{fifteenh}  
for the generalized cosecant numbers of the form, $c_{2v,v-\ell}$, where $\ell=2,3, \dots$. These have yet to be developed.

Finally, it should be mentioned that the series on the lhs of \eqref{sixteen} to \eqref{twenty} also represent
specific values of the generalized harmonic numbers, which are defined as $H_{n,r}=\sum_{k=1}^{n}1/k^r$ \cite{Sond2017}. 
In particular, for the case of $r=2$ given by \eqref{sixteen} the numbers are known as Wolstenholme numbers, which appear 
as Sequences A007406, A007408, A11354 and A123751 in the online encyclopedia of integer sequences \cite{Sl2017}.

\section{Acknowledgement}
The author is grateful to Professor Carlos M. da Fonseca, Kuwait University, without whose support and encouragement this
paper would not have proceeded.

\end{document}